\documentclass{amsart}
\usepackage{a4wide,amssymb,amsrefs,hyperref}

\newcommand{\mnfld}{\mathcal{M}}
\newcommand{\Z}{\mathbb{Z}}
\newcommand{\R}{\mathbb{R}}
\newcommand{\HHH}{\mathbb{H}}
\newcommand{\SL}{\operatorname{SL}}
\newcommand{\spc}{\mathfrak{E}}
\newcommand{\spcc}{\mathfrak{F}}
\newcommand{\cnst}{C}
\newcommand{\unq}{\operatorname{unq}}
\newcommand{\stds}{\mathcal{P}}
\newcommand{\rest}{\lvert}
\newcommand{\funct}{\mathfrak{H}}      
\newcommand{\functb}{\mathfrak{F}}        
\newcommand{\umd}{\functb_{\operatorname{umg}}}
\newcommand{\Weyl}{W}
\newcommand{\cusp}{\operatorname{cusp}}
\newcommand{\Siegel}{\mathcal{S}}
\newcommand{\autspace}{\mathcal{X}}
\newcommand{\zspace}{\mathcal{Z}}

\newcommand{\modulus}{\delta}
\newcommand{\srts}{\Delta}
\newcommand{\C}{\mathbb{C}}

\newcommand{\A}{\mathbb{A}}
\newcommand{\z}{\mathfrak{z}}
\newcommand{\aaa}{\mathfrak{a}}
\newcommand{\Ht}{H}
\newcommand{\AF}{\mathcal{A}}
\newcommand{\iii}{{\rm i}}
\newcommand{\bs}{\backslash}
\newcommand{\sprod}[2]{\left\langle#1,#2\right\rangle}
\newcommand{\abs}[1]{\left|{#1}\right|}
\newcommand{\norm}[1]{\lVert#1\rVert}

\newcommand{\Exp}{\mathcal{E}}
\newcommand{\Img}{\operatorname{Im}}
\newcommand{\Rel}{\operatorname{Re}}
\newcommand{\GL}{\operatorname{GL}}
\newcommand{\Ind}{\operatorname{Ind}}
\newcommand{\disc}{\operatorname{disc}}
\newcommand{\Polyexp}{\mathcal{PE}}
\newcommand{\unit}{\operatorname{unit}}
\newcommand{\temp}{\operatorname{temp}}
\newcommand{\weak}{\operatorname{weak}}
\newcommand{\inner}[2]{\big\langle#1,#2\big\rangle}
\newcommand{\swrz}{\mathcal{S}}
\newcommand{\opr}{\partial}
\newcommand{\posrts}[2]{D_{#1,+}(#2)}
\newcommand{\newchar}{\chi}
\newcommand{\SO}{\operatorname{SO}}

\newtheorem{theorem}{Theorem}[section]

\newtheorem{proposition}[theorem]{Proposition}

\theoremstyle{definition}

\theoremstyle{remark}
\newtheorem{remark}[theorem]{Remark}

\numberwithin{equation}{section}

\begin{document}

\title[Eisenstein series]{Some perspectives on Eisenstein series}
\author{Erez Lapid}
\address{Department of Mathematics, Weizmann Institute of Science, Rehovot 7610001, Israel}
\email{erez.m.lapid@gmail.com}


\date{\today}

\begin{abstract}
We review some topics in the analytic theory of Eisenstein series, including meromorphic continuation, $L^2$-spectral expansion
and Fourier coefficients. We also discuss some open problems.
\end{abstract}

\maketitle

\setcounter{tocdepth}{1}
\tableofcontents

\section{Introduction}
The goal of these notes is to briefly describe some developments in the theory of Eisenstein series
in the last 40 years or so. We will also discuss some open questions and problems.

The role of Eisenstein series in the theory of automorphic forms is ubiquitous\footnote{As of April 1, 2022
there are 855 articles on the MathSciNet database which have ``Eisenstein series'' in their title.}
and it would be
fatuous (and certainly beyond the author's ability) to try to provide a comprehensive account.
After all, any serious use of spectral theory, whether in the trace formula or other applications,
inevitably involves Eisenstein series.
Instead, we focus on a few select topics in the analytic theory, their choice heavily biased by personal taste and experience.

We start by sketching Bernstein's proof of the meromorphic continuation of Eisenstein series,
which simplifies the earlier proofs of Selberg and Langlands.

This will be followed by an outline of a proof of the spectral decomposition of $L^2(G(F)\bs G(\A))$,
consolidating ideas of Bernstein, Delorme and the author, once again simplifying Langlands's proof.

Needless to say, both proofs rely on basic results of Selberg and Langlands,
as well as Harish-Chandra and Gelfand--Piatetski-Shapiro.

We will then review, rather telegraphically, some constructions involving Eisenstein series, including
the Langlands--Shahidi method,
the descent construction of Ginzburg--Rallis--Soudry and the double descent of Ginzburg--Soudry.

We will end with a (rather amorphous) list of open problems.

Many other topics will be omitted altogether, including arithmetic applications, function field and geometric analogues.
In particular, we will not discuss at all congruences between Eisenstein series and cusp forms,
a vast subject which started with the pioneering works of Swinnerton-Dyer, Serre, Ribet, Katz, Mazur, Wiles and others \cites{MR0406931, MR0404145, MR419403, MR0506271, MR488287, MR742853}.

The prototype of an Eisenstein series is
\[
E(z;s)=\tfrac12\sum_{(m,n)\in\Z^2\mid\gcd(m,n)=1}\tfrac{(\Img z)^{s+\frac12}}{\abs{mz+n}^{2s+1}},\ \ \Img z>0,
\]
introduced by Maass in his seminal paper \cite{MR31519}, where the automorphic forms that bear his name also originated.
The series converges absolutely for $\Rel s>\frac12$.
A more flexible definition, which is part of the modern theory of Eisenstein series, gives (for $s=k-\frac12$) the elliptic modular forms
of weight $2k$ that were originally considered by Eisenstein.\footnote{We refer to Andre Weil's book review of Eisenstein's collected works
(Bulletin of the American Mathematical Society, 1976) for a lucid account on Eisenstein's life and work and to \cite{MR1723749} for a
more in-depth discussion.}

The meromorphic continuation for $s\in\C$ was achieved in part by Roelcke \cite{MR82562} and in full
by Selberg \cite{MR0088511} who recognized their importance to spectral theory.
Langlands generalized Selberg's results to the higher rank case and gave a spectral decomposition
of the $L^2$-space of an automorphic quotient in terms of Eisenstein series built from automorphic forms
belonging to the discrete spectrum of automorphic quotients pertaining to parabolic subgroups \cites{MR0579181, MR1361168}.

The theory of Eisenstein series is not limited to arithmetic quotients.\footnote{Eisenstein series for non-arithmetic
lattices of $\SL_2(\R)$ are important for instance for dynamics on Veech surfaces, including regular polygons. See \cites{MR1005006,MR4254262}.}
However, we will restrict ourselves to the adelic setup, i.e., to automorphic quotients by congruence (arithmetic) subgroups.
Also, we will only consider the number field case.

\section{Meromorphic continuation of Eisenstein series}
\subsection{}
Let $G$ be a reductive group over a number field $F$ and let $P=M\ltimes U$ be a parabolic subgroup
defined over $F$.
The Eisenstein series in its general form is built from an automorphic form on $M(F)\bs M(\A)$.
They depend on a spectral parameter, namely a character of $M(\A)$, trivial on $M(F)$.
These characters form a complex analytic space whose connected components are linear spaces of constant dimension,
namely the rank of the lattice of $F$-rational characters of $M$.
The Eisenstein series converges and is holomorphic when the real part lies sufficiently deep
in the positive Weyl chamber, determined by $P$.
The first order of business is to meromorphically continue the Eisenstein series
and to provide functional equations for them. The functional equations are in terms of intertwining operators,
which satisfy functional equations of their own.

Langlands's approach to the meromorphic continuation of Eisenstein series
was to start with the case where the inducing automorphic form is cuspidal (on the Levi subgroup).
This case was handled more or less by Selberg's ideas in the rank one case.
(Subsequently, Selberg honed his method of proof and gave some hints in \cite{MR0176097}.)
In one of his greatest achievements, Langlands showed that the discrete spectrum is obtained from the
cuspidal spectrum of Levi subgroups by means of residues of Eisenstein series.
He then deduced the meromorphic continuation of Eisenstein series in case the inducing data
is square-integrable and completed the spectral decomposition of $L^2(G(F)\bs G(\A))$.
Subsequently, Franke showed (using the spectral decomposition) that every automorphic form can be written as
a linear combination of Laurent coefficients of Eisenstein series induced from cuspidal automorphic forms \cite{MR1603257}.
This implies the meromorphic continuation of Eisenstein series induced from an arbitrary automorphic form.

In the 1980s Joseph Bernstein came up with a new and simple idea to meromorphically continue the Eisenstein series.
His method works for any automorphic form, not necessarily cuspidal or discrete.
Arguably, Bernstein's main insight was that the meromorphic
continuation of Eisenstein series is completely separate from spectral theory.
In hindsight, the ostensibly inextricable bond between Eisenstein series and $L^2$ theory
(``the bed of Procrustes'' as Langlands refers to it) was a major source of difficulty in the proofs
of Selberg and Langlands \cites{MR0249539, MR993313}.

Let me outline Bernstein's approach, referring to \cite{1911.02342} for more details.

\subsection{}
Bernstein's proof is based on the ``principle of meromorphic continuation'' (PMC),
pertaining to holomorphic systems of linear equations (HSLE).
By that we mean a family $\Xi(s)$, where $s$ ranges over a connected complex analytic manifold $\mnfld$, of linear equations
\[
\sprod{\lambda_i(s)}{v}=c_i(s),\ \ i\in I
\]
on $v$ in a Hausdorff, locally convex topological vector space (LCTVS) $\spc$ over $\C$, whose coefficients depend holomorphically on $s$. That is,
for every $i\in I$
\begin{itemize}
\item $c_i(s)$ is a holomorphic function on $\mnfld$.
\item $\lambda_i(s)$, $s\in\mnfld$ is a family of continuous linear functionals on $\spc$
such that for every $u\in\spc$ the function $\sprod{\lambda_i(s)}{u}$ is holomorphic in $s$.
\end{itemize}
In real life, a more compact and natural way to write a HSLE is
\[
A_i(s)v=u_i(s)
\]
where for each $i$, $s\in\mnfld\mapsto A_i(s)$ is a holomorphic family of continuous operators from $\spc$ to a LCTVS $\spcc_i$
and $u_i(s)$ is a holomorphic family of vectors in $\spcc_i$.

We assume that the set $\mnfld_{\unq}$, where $\Xi(s)$ admits a unique solution $v(s)$, has a nonempty interior.
The PMC states that under suitable additional conditions, $\mnfld_{\unq}$ is the complement of a closed analytic
subset of $\mnfld$ and the function $v(s)$ is holomorphic in $\mnfld_{\unq}$ and extends to a meromorphic function on $\mnfld$.
In other words, locally in $\mnfld$, there exists a holomorphic function $f$ with no zeroes on $\mnfld_{\unq}$,
such that $f(s)\sprod{\mu}{v(s)}$ is holomorphic for every $\mu$ in the continuous dual of $\spc$.

The easiest case where PMC holds is when $\Xi(s)$ is locally of finite type (LFT). This means that locally,
the set of solutions of $\Xi(s)$ is contained in the image of an injective linear operator
from a finite-dimensional vector space $L$ to $\spc$ that depends analytically on $s$.
Indeed, this is a simple application of Cramer's rule.
In turn, the basic tool for proving local finiteness is Fredholm theory.
Suppose for instance that $\mu_s,\nu_s:\spc\rightarrow\spcc$, $s\in\mnfld$ are two analytic families
of bounded operators between Hilbert spaces.
Suppose that for all $s\in \mnfld$, $\mu_s$ is a strict embedding and $\nu_s$ is a compact operator.
Then, the homogeneous equation
\[
\mu_sv=\nu_sv
\]
on $v\in\spc$ is LFT. Such a system will be called ``of Fredholm type''.
A practical way to prove that an operator is compact is to show that it is Hilbert--Schmidt.
For instance, a bounded operator $A:\spc\rightarrow L^2(X,\mu)$ from a Hilbert space $\spc$ to an $L^2$-space
is Hilbert--Schmidt if and only if for almost all $x\in X$ the evaluation $ev_x(u)=Au(x)$ is bounded on $\spc$
and $x\mapsto\norm{ev_x}\in L^2(X,\mu)$.

\subsection{}
Next, we describe the HSLE that is used in the proof of
meromorphic continuation of Eisenstein series,
starting with the $\SL_2$ case, which in fact contains much of the ideas for the general case.

Let $G=\SL_2(\R)$ and let $\HHH=G/K$ be the upper half-plane (where $K=\SO(2)$) with the left $G$-action by M\"obius transformations.
Let $\autspace=\Gamma\bs\HHH$ where $\Gamma=\SL_2(\Z)$.
Consider the Eisenstein series
\begin{equation} \label{eq: ezs}
E(z;s)=\sum_{\gamma\in\pm\Gamma_\infty\bs\Gamma}\Img(\gamma z)^{s+\frac12}=
\tfrac12\sum_{(m,n)\in\Z^2\mid\gcd(m,n)=1}\tfrac{y^{s+\frac12}}{\abs{mz+n}^{2s+1}}
\end{equation}
where $\Gamma_\infty=\{\left(\begin{smallmatrix}1&n\\0&1\end{smallmatrix}\right)\mid n\in\Z\}\subset\Gamma$.
The series converges absolutely for $\Rel s>\frac12$ and defines a function of uniform moderate growth on $\autspace$.
(Denote this space by $\umd(\autspace)$ -- cf.\ \cite{MR1482800} for some basic facts.)

The HSLE $\Xi(s)$, $s\in\C$ on $\psi\in\umd(\autspace)$ is given by the following three sets of linear equations.
\begin{subequations}
\begin{gather}
\label{eq: 1cls} \delta(h)\psi=\hat h(s)\psi\ \ \ \forall h\in C_c^\infty(G//K),\\
\label{eq: 2csls} (T_a-a^{-s})(\cnst\psi(y)-y^{s+\frac12})\equiv0\ \ \forall a>0,\\
\label{eq: 3cls} (\psi,f)_{\autspace}=0\ \ \text{for every cusp form $f$ on $\autspace$}.
\end{gather}
\end{subequations}

The notation is the following.
\begin{itemize}
\item $\delta(h)$ denotes the right action (on functions on $\HHH$) by bi-$K$-invariant functions on $G$.
\item $\hat h(s)$ is the eigenvalue of the eigenfunction $(\Img z)^{s+\frac12}$ under $\delta(h)$.\footnote{$\hat h(s)$
can be computed explicitly. We only need to know that it is entire and that for every $s\in\C$
there exists $h\in C_c^\infty(G//K)$ such that $\hat h(s)\ne0$.}
\item $\cnst f$ denotes the constant term
$\cnst f(y)=\int_{\Z\bs\R}f(x+\iii y)\ dx$, $y>0$.
\item $T_a$ is the normalized shift operator
$T_af(y)=a^{-\frac12}f(ay)$ on functions on $\R_{>0}$.
\item $(\cdot,\cdot)_\autspace$ is the pairing with respect to the measure $\mu=\frac{dx\, dy}{y^2}$ on $\autspace$.
\end{itemize}

The fact that the Eisenstein series satisfies $\Xi(s)$ for $\Rel s>\frac12$ is straightforward.
Note that the non-homogeneous equation \eqref{eq: 2csls} amounts to the equation
\[
\cnst E(y;s)=y^{s+\frac12}+m(s)y^{-s+\frac12},\ \ \ y>0
\]
for some function $m(s)$. However, a priori we do not need to know \emph{anything} about $m(s)$.
In fact, $m(s)$ can be expressed in terms of the Riemann zeta function $\zeta(s)$, and ultimately the proof
provides an alternative approach for the meromorphic continuation of $\zeta(s)$
(albeit with no information about its poles).

In order to apply PMC, we show the following two statements.
\begin{enumerate}
\item For $\Rel s>\frac12$, the Eisenstein series $E(z;s)$ is the unique function of uniform moderate growth
satisfying the equations \eqref{eq: 2csls} and \eqref{eq: 3cls}.
\item The HSLE \eqref{eq: 1cls} and \eqref{eq: 2csls} (already for a single $a>0$) is LFT.
\end{enumerate}

For the first part, the point is that if $\psi'\in\umd(\autspace)$ and $\cnst\psi'$ is proportional to $y^{\frac12-s}$,
then $\psi'$ is bounded for $y\ge\frac12$ (say) since $\psi'-\cnst\psi'$ is rapidly decreasing.
Hence, $\psi'$ is bounded on $\autspace$ since it is $\Gamma$-invariant.
Therefore, $\cnst\psi'$ is bounded on $\R_{>0}$. By the condition on $\psi'$ and $s$, we conclude that $\cnst\psi'\equiv0$, so that
$\psi'$ is cuspidal. Together with \eqref{eq: 3cls} this implies that $\psi'\equiv0$.

In order to show the second part, we pass to an auxiliary HSLE on a Hilbert space and replace the complicated space
$\autspace$ by a simpler one which approximates it at the cusp.

More precisely, let $\zspace=\Gamma_\infty\bs\HHH$ and let $\Siegel$ be the truncated space $y>c_0$
where $c_0>0$ is chosen so that the projection $p:\Siegel\rightarrow\autspace$ is onto.
For $N>0$ consider the Hilbert space $\funct^N(\Siegel)=L^2(\Siegel;y^{-2N}\mu)$.
The pullback $f\mapsto\tilde f$ by $p$
gives rise to a Hilbert space $\funct^N(\autspace)$ of functions on $\autspace$ which strictly embeds in $\funct^N(\Siegel)$
and whose topology is independent of the choice of $c_0$. Moreover, $\umd(\autspace)$ is the union over $N$
of the smooth part of $\funct^N(\autspace)$.

Any $f\in\funct^N(\Siegel)$ admits an orthogonal decomposition
$f=\cnst f+f_{\cusp}$ where $\cnst f_{\cusp}\equiv0$.

Fix $a>1$ and $h\in C_c^\infty(G//K)$.
We claim that the following HSLE $\tilde\Xi^N(s)$ on $f\in\funct^N(\autspace)$
is of Fredholm type for $\abs{\Rel s}<N$ and $\hat h(s)\ne0$.
\begin{align*}
\hat h(s)\tilde f_{\cusp}&=\widetilde{\delta(h)f}_{\cusp},\\
\hat h(s)(\cnst f)\rest_{[c_0,c_0a^2]}&=\cnst(\delta(h)f)\rest_{[c_0,c_0a^2]},\\
(T_a-a^s)(T_a-a^{-s})(\cnst f)&=0.
\end{align*}
More precisely, the operator
\begin{align*}
\funct^N(\autspace)&\rightarrow\funct^N(\Siegel)\oplus L^2(\R_{>c_0},y^{-2N}\tfrac{dy}{y^2})\oplus L^2([c_0,c_0a^2])\\
f&\mapsto (\tilde f_{\cusp},(T_a-a^s)(T_a-a^{-s})(\cnst f),(\cnst f)\rest_{[c_0,c_0a^2]})\\
\text{(resp., }f&\mapsto(\widetilde{\delta(h)f}_{\cusp},0,\cnst(\delta(h)f)\rest_{[c_0,c_0a^2]})\text{ )}
\end{align*}
is a strict embedding (resp., compact (and in fact Hilbert--Schmidt)).

The first statement boils down to the elementary fact that the operator
\begin{align*}
L^2(\R_+,e^{-2rx}\ dx)&\rightarrow L^2(\R_+,e^{-2rx}\ dx)\oplus L^2([0,1]),\\
f&\mapsto (f(x+1)-\lambda f(x),f\rest_{[0,1]})
\end{align*}
is a strict embedding provided that $e^r>\abs{\lambda}$.
The second statement is standard.

It follows that $\tilde\Xi^N(s)$ is LFT in the domain above.
The local finiteness of \eqref{eq: 1cls}+\eqref{eq: 2csls} is deduced by an approximate identity argument since for any $h\in C_c^\infty(G//K)$,
$\delta(h)$ defines a continuous operator from $\funct^N(\autspace)$ to $\umd(\autspace)$.

\subsection{}
Let us describe in broad strokes the HSLE in the general case.

Fix a minimal parabolic subgroup $P_0$ of $G$ defined over $F$ and denote by $\stds$ the finite set of standard
parabolic subgroups of $G$ defined over $F$.

Fix $P=M\ltimes U\in\stds$ and an automorphic form $\varphi$ on $\autspace_P=P(F)U(\A)\bs G(\A)$.
Denote by $m_P(g)$ the $M(\A)$-part of $g\in G(\A)$ in the Iwasawa decomposition with respect to a suitable (fixed)
choice of a maximal compact subgroup of $G(\A)$.

Let $\srts_P'\subset\srts_0$ be the set of simple roots whose restriction to $\operatorname{Lie}(U)$ is nontrivial.

For a quasi-character $\lambda$ of $M(\A)/M(\A)^1$, the Eisenstein series
\[
E(g,\varphi,\lambda)=\sum_{\gamma\in P(F)\bs G(F)}\varphi(\gamma g)m_P(\gamma g)^\lambda
\]
converges absolutely if $\Rel\lambda\gg0$ (that is, $\Rel\sprod{\lambda}{\alpha^\vee}\gg0$ for all $\alpha\in\srts_P'$)
and defines an automorphic form on $\autspace_G=G(F)\bs G(\A)$.

Denote by $\Weyl$ the Weyl group of $G$.
Then, for any $Q=L\ltimes V\in\stds$, we have
\[
\Exp_Q^{\cusp}(E(\varphi,\lambda))\subset
\{w(\lambda+\mu)\mid w\in\Weyl^{\supset Q}(P), \mu\in\Exp_{P_w}^{\cusp}(\varphi)\}\text{       where}
\]
\begin{itemize}
\item $\Exp_R^{\cusp}(\phi)$ is the set (with multiplicities) of cuspidal exponents of $\phi$ along a parabolic subgroup $R\subset P$.
\item $\Weyl^{\supset Q}(P)=\{w\in\Weyl\text{ right $\Weyl_M$-reduced}\mid wMw^{-1}\supset L\}$.
\item For any $w\in \Weyl^{\supset Q}(P)$, $P_w$ is the standard parabolic subgroup of $P$ with Levi subgroup $w^{-1}Lw$.
\end{itemize}
This follows from the computation of the constant term of $E(\varphi,\lambda)$ along $Q$ as a sum of contributions
over the double cosets
$P(F)\bs G(F)/Q(F)$, or equivalently over $w\in\Weyl_L\bs\Weyl/\Weyl_M$ using the Bruhat decomposition.
(It is analogous to the computation of the Jacquet module of an induced representation in the local case, i.e.,
the geometric lemma of Bernstein--Zelevinsky.)
While most of the terms involve nontrivial intertwining operators and should be treated as unknown
(apart from their exponents),
the term corresponding to $w=e$ is simply the constant term $\cnst_Q\varphi$ of $\varphi$ itself (interpreted as $0$ if $Q\not\subset P$).
Thus, denoting by $\cnst_Q^{\cusp}\phi$ the cuspidal projection of $\cnst_Q\phi$,
the exponents of the difference $\cnst_Q^{\cusp}E(\varphi,\lambda)-\cnst_Q^{\cusp}\varphi$ are contained (as a multiset) in
\[
\{w(\lambda+\mu)\mid w\in\Weyl^{\supset Q}(P)\setminus\{e\}, \mu\in\Exp_{P_w}^{\cusp}(\varphi)\}.
\]
Moreover, this property determines $E(\varphi,\lambda)$ uniquely, at least if $\Rel\lambda\gg0$.
This is a consequence of the following assertion which is a strengthening of a basic result of Langlands, and is interesting in its own right.
\begin{quote}
Let $\phi$ be an automorphic form, not identically zero.
Then, there exists $Q\in\stds$ and $\lambda\in\Exp_Q^{\cusp}(\phi)$
such that $\Rel\lambda+\rho_Q$ lies in the closure of the positive Weyl chamber.
\end{quote}
This is proved by reducing it to a corank one statement, where the argument is very similar to the $\SL_2$ case
explained above.

The deduction of the uniqueness statement for $E(\varphi,\lambda)$ follows by observing that if $\Rel\lambda\gg0$
and $w\in\Weyl^{\supset Q}(P)\setminus\{e\}$, then $w\Rel\lambda$ is far from the positive Weyl chamber of $Q$.

For the local finiteness we need an additional set of equations. They are of the form
\[
\delta(h_i(\lambda))\psi=c_i(\lambda)\psi,\ \ i\in I
\]
where for each $i\in I$, $h_i(\lambda)$ is a holomorphic family of functions in $C_c^\infty(G(\A))$
and $c_i(\lambda)$, $i\in I$ are holomorphic functions with no common zeros.

By a basic result of Harish-Chandra, such equations are satisfied by any holomorphic family of automorphic forms.

To show local finiteness, assume for simplicity that the center of $G$ is $F$-anisotropic.
We show that the following HSLE on $f\in\umd(\autspace)$ is LFT:
\begin{equation} \label{eq: HSLE1}
\delta(h(s))f=f,\ \ D_{\alpha}(s)(T_{a_\alpha})(\cnst_{P_{\alpha}}f)=0\text{ for every }\alpha\in\srts_0.
\end{equation}
Here, $h(s)$ is a holomorphic family of functions in $C_c^\infty(G(\A))$ (for $s$ in a complex analytic manifold) and for every $\alpha\in\srts_0$
\begin{itemize}
\item $D_\alpha(s)$ is a holomorphic family of monic polynomials in one variable of degree $m_\alpha$.
\item $P_\alpha$ is the maximal parabolic subgroup of $G$ corresponding to $\alpha$.
\item $T_{a_\alpha}$ is the normalized left translation by a fixed element $a_\alpha$ of the center of the Levi part
    of $P_\alpha$ such that $\abs{\alpha(a_\alpha)}>1$.
\end{itemize}

As in the $\SL_2$ case, in order to prove this, we pass to an auxiliary HSLE $\Xi^\lambda(s)$, described below,
on a certain Hilbert space $\funct^\lambda(\autspace)$. Here $\lambda\in\aaa_0^*$ is a parameter.

Fix a Siegel domain $\Siegel$ of $P_0(F)\bs G(\A)$ such that the projection $p:\Siegel\rightarrow\autspace$ is onto.
Define a weighted $L^2$-space
\[
\funct^\lambda(\Siegel)=L^2(\Siegel,m_{P_0}(x)^{-2\lambda}\ dx).
\]
The pullback $f\mapsto f^\Siegel$ by $p$ gives rise to a Hilbert space $\funct^\lambda(\autspace)$ of functions on $\autspace$ which strictly embeds in $\funct^\lambda(\Siegel)$ and whose topology is independent of the choice of $\Siegel$. This space was considered by Franke. We have
$\umd(\autspace)=\cup_{\lambda}\funct^\lambda_\infty(\autspace)$.

The space $\funct^\lambda(\Siegel)$ admits a Harish-Chandra decomposition \cite{MR0232893}
\[
\funct^\lambda(\Siegel)=\oplus_{P\in\stds}\funct^\lambda_{\cusp}(\Siegel_P)
\]
where $\Siegel_P$ is the image of $\Siegel$ under the projection to $U(\A)P_0(F)\bs G(\A)$.
For any $f\in\funct^\lambda(\autspace)$ denote by $f_P^{\Siegel}$, $P\in\stds$ the components of $f^{\Siegel}$ with
respect to this decomposition.

The HSLE $\Xi^\lambda(s)$ on $f\in \funct^\lambda(\autspace)$ consists of the equations
\[
D_\alpha(s)(T_{a_\alpha})(f_P^{\Siegel})=0\text{ for all }\alpha\in\srts_P',\ \
f_P^{\Siegel}\rest_{\Siegel_P'}=(\delta(h(s))f)_P^{\Siegel}\rest_{\Siegel_P'}
\]
for every $P\in\stds$, where $\Siegel_P'=\Siegel_P\setminus\cup_{\alpha\in\srts_P'}a_\alpha^{m_\alpha}\Siegel_P$.

It suffices to show that $\Xi^\lambda(s)$ is LFT provided that $\lambda$ is sufficiently positive.
In fact, it is of Fredholm type. Namely, as in the $\SL_2$ case, it is easy to see that for every $P$, the operator
\[
\funct^\lambda(\Siegel_P)\rightarrow \funct^\lambda(\Siegel_P)^{\srts_P'}\oplus\funct^\lambda(\Siegel_P'),\ \
f\mapsto \big((D_{\alpha}(s)(T_{a_\alpha})(f))_{\alpha\in\srts_P'},f\rest_{\Siegel_P'}\big)
\]
is a strict embedding provided that $e^{\sprod{\lambda}{\Ht_0(a_\alpha)}}>\abs{r}$ for every root $r$ of $D_\alpha(s)$
(which also implies that the solutions of \eqref{eq: HSLE1} are in $\funct^\lambda(\autspace)$).
On the other hand, by a variant of a familiar result of Gelfand and Piatetski-Shapiro \cite{MR0159899}, the operator
\[
\funct^\lambda(\autspace)\rightarrow\funct^\lambda(\Siegel_P'),\ \ f\mapsto(\delta(h(s))f)_P^{\Siegel}\rest_{\Siegel_P'}
\]
is compact (and in fact, Hilbert--Schmidt).

\subsection{}
A natural question, for which we do not have an answer, is whether the method of proof can be extended to deal with Eisenstein series induced
from smooth automorphic forms (that is, $\z$-finite but not necessarily $K$-finite smooth functions of uniform moderate growth).
Such functions are no longer eigenfunctions under convolution by a compactly supported smooth function on $G(\A)$.
Nonetheless, one can optimistically hope that the PMC is still applicable
for a suitable system, even though it would not be LFT (let alone of Fredholm type) in the strict sense.
At the very least, this would entail formulating and proving a more flexible version of the PMC.
We mention that there are other approaches to extend the meromorphic continuation from the $K$-finite case to the smooth case \cites{MR2402686, MR3219530, MR4389792}.

Another open problem, already in the $K$-finite case, is whether the PMC can establish, at least in principle,
that the solution is a meromorphic function
of finite order. Once again, in the particular case of Eisenstein series, this is known by other methods \cite{MR1025165}.

\section{Spectral decomposition of \texorpdfstring{$L^2(G(F)\bs G(\A))$}{L2}}

Once the meromorphic continuation of Eisenstein series is established, it is not difficult
to obtain the spectral decomposition of $L^2(G(F)\bs G(\A))$.
We will present the argument following Bernstein \cite{MR1075727}, Delorme \cite{2006.12893} and the author \cites{MR2767521, MR3001807}.
It is analogous to the local case \cites{MR1989693, MR1626757} and provides a simplification of the original argument of Langlands.

\subsection{}
For any $P\in\stds$ denote by $\AF_P$ the space of automorphic forms on $\autspace_P$.
Let
\[
\AF_P^2=\{\varphi\in\AF_P\mid\modulus_P^{-\frac12}\varphi(\cdot g)\in
L^2(A_MM(F)\bs M(\A))\ \forall g\in G(\A)\}.
\]
This is an inner product space with respect to integration over $A_M\bs\autspace_P$.
Its completion is
\[
\Ind_{P(\A)}^{G(\A)}L^2_{\disc}(A_MM(F)\bs M(\A)).
\]

By the functional equations, for any $w\in\Weyl(P,Q)$ the intertwining operators
\[
M(w,\lambda):\AF_P^2\rightarrow\AF_Q^2,\ \ \lambda\in\aaa_{P,\C}^*
\]
are unitary (and in particular, holomorphic) on $\iii\aaa_P^*$ and extend to isometries
\[
\Ind_{P(\A)}^{G(\A)}L^2_{\disc}(A_MM(F)\bs M(\A))\rightarrow
\Ind_{Q(\A)}^{G(\A)}L^2_{\disc}(A_LL(F)\bs L(\A)).
\]

\begin{theorem}[Langlands] \label{thm: Langlands}
The bilinear map
\[
f,\varphi\in C_c^\infty(\iii\aaa_P^*)\times\AF_P^2\mapsto
E_{f\otimes\varphi}=\int_{\iii\aaa_P^*}f(\lambda)E(\varphi,\lambda)\ d\lambda
\]
(suitably normalized) extends to an isometry of Hilbert spaces
\[
\Theta:\big(\bigoplus_PL^2(\iii\aaa_P^*;\Ind_{P(\A)}^{G(\A)}L^2_{\disc}(A_MM(F)\bs M(\A))\big)^{\Weyl}\rightarrow
L^2(\autspace_G)
\]
where the space on the left-hand side is the closed subspace defined by the relations\footnote{It is better to view
it as the coinvariants under a groupoid.}
\[
f^Q(w\lambda)=M(w,\lambda)f^P(\lambda)\ \ \lambda\in\iii\aaa_P^*,\ w\in\Weyl(P,Q).
\]
\end{theorem}

\subsection{}
We say that an automorphic form $\varphi\in\AF_P$ is \emph{tempered} if there exists $k$ such that
\[
\abs{\varphi(g)}\ll e^{\sprod{\rho}{H(g)}}(1+\norm{H(g)})^k
\]
for all $g\in A_M\Siegel$. Equivalently, every exponent $\lambda\in\Exp_Q(\varphi)$, $Q\subset P$ is subunitary
(i.e., $\Rel\lambda$ is a sum of simple roots $\srts_Q^P$ with non-positive coefficients.)
We denote by $\AF_P^{\temp}$ the space of tempered automorphic forms.
Note that $\AF_P^2\subset\AF_P^{\temp}$.

In particular, if $\varphi\in\AF_P^2$ then $E(\varphi,\lambda)\in\AF_G^{\temp}$ for all $\lambda\in\iii\aaa_P^*$ at which
$E(\varphi,\lambda)$ is holomorphic.
(Eventually, the last condition is redundant.)

For $\varphi\in\AF_G^{\temp}$, define the \emph{weak constant term} $\varphi_P^{\weak}$ to be the part of $\varphi_P$ corresponding to the exponents $\lambda$ with $\Rel\lambda=0$. Then $\varphi_P^{\weak}\in\AF_P^{\temp}$. Note that if $\varphi\in\AF_G^2$, then $\varphi_P^{\weak}\equiv0$ for all proper $P$.

For any $\phi\in\AF_P^2$ and $\psi\in\AF_P^{\temp}$ the tempered distribution $p(\phi,\psi)$ on $\aaa_P$
\[
f\in\swrz(\aaa_P)\mapsto\inner{\phi\cdot(f\circ\Ht_P)}{\psi}_{\autspace_P}
\]
is a polynomial exponential whose exponents are unitary and contained in the set $\Exp_P(\varphi)+\overline{\Exp_P(\psi)}$.

For any $p\in\Polyexp^{\unit}(\aaa_P)$ (i.e., a polynomial exponential with unitary exponents) denote by $\opr(p)$ the convolution
(on $C_c^\infty(\iii\aaa_P^*)$) by the Fourier transform of $p$, viewed as a finitely supported distribution on $\iii\aaa_P^*$.

It is a basic fact that for any $\psi\in\AF_G^{\temp}$,
\[
\psi\equiv0\iff p(\varphi,\psi_P^{\weak})\equiv0\text{ for every $P\in\stds$ and }\varphi\in\AF_P^2.
\]
In view of Bernstein's result on the support of the Plancherel measure \cite{MR1075727},
Theorem \ref{thm: Langlands} is essentially equivalent to the following.
\begin{theorem}\cite{2006.12893} \label{thm: main2}
Let $\varphi\in\AF_P^2$. Then,
\begin{enumerate}
\item For any $\psi\in\AF_G^{\temp}$ and $f\in C_c^\infty(\iii\aaa_P^*)$ we have
\[
\inner{\int_{\iii\aaa_P^*}f(\lambda) E(\varphi,\lambda)\ d\lambda}{\psi}_{\autspace_G}=
\opr(p(\varphi,\psi_P^{\weak}))f(0).
\]
\item Assume that $\psi=E(\varphi',\mu)$ where $\varphi'\in\AF_Q^2$ and $\mu\in\iii\aaa_Q^*$.
Then,
\[
p(\varphi,\psi_P^{\weak})=\sum_{w\in\Weyl(Q,P)}p(\varphi,(M(w,\mu)\varphi')_{w\mu}).
\]
Thus,
\[
\inner{\int_{\iii\aaa_P^*}f(\lambda) E(\varphi,\lambda)\ d\lambda}{\psi}_{\autspace_G}=
\sum_{w\in\Weyl(Q,P)}f(w\mu)\inner{\varphi}{M(w,\mu)\varphi'}_{A_M\bs\autspace_P}.
\]
\end{enumerate}
\end{theorem}

We note that the left-hand side in the first part is well defined since the ``wave packet''
$\int_{\iii\aaa_P^*}f(\lambda) E(\varphi,\lambda)\ d\lambda$ is in the Harish-Chandra Schwartz space of $\autspace_G$.

\subsection{}

One way to prove Theorem \ref{thm: main2} is using Arthur's truncation operator $\Lambda^T$ \cite{MR558260}.
(In \cite{2006.12893} a more straightforward truncation is used.)
Assume for simplicity that the center of $G$ is $F$-anisotropic.
A basic fact is that for any $\phi,\psi\in\AF_G$ the function
\[
\inner{\Lambda^T\phi}{\psi}_{\autspace_G}=\inner{\phi}{\Lambda^T\psi}_{\autspace_G}
\]
is a polynomial exponential in $T$ whose exponents are contained in the set
\[
\cup_P(\Exp_P(\varphi)+\overline{\Exp_P(\psi)}).
\]
In particular, if $\phi,\psi\in\AF_G^{\temp}$, then the exponents are subunitary.
In this case, we denote by $\inner{\Lambda^T\phi}{\psi}_{\autspace_G}^{\unit}$
the unitary part of this polynomial exponential.

For now, we pretend that we know that the Eisenstein series $E(\varphi,\lambda)$ is holomorphic for $\lambda\in\iii\aaa_P^*$
for every $\varphi\in\AF_P^2$. (See Remark \ref{rem: eisholom} below.)

\begin{proposition}(cf.\ \cite{2006.12893}) \label{prop: maincor}
Let $\psi\in\AF_G^{\temp}$ and $\varphi\in\AF_P^2$. Fix $T\in\aaa_0$.
Then, we have an equality of distributions on $\iii\aaa_P^*$
\begin{equation} \label{eq: main34}
\inner{\Lambda^TE(\varphi,\lambda)}{\psi}_{\autspace_G}^{\unit}=
\sum_Q\sum_{w\in\Weyl(P,Q)}\opr(p(M(w,w^{-1}\cdot)\varphi,\psi_Q^{\weak}))(\widehat{\chi_w^T})\circ w
\end{equation}
where $\chi_w^T$ (viewed as a tempered distribution on $\aaa_Q$) is
$(-1)^{\#\{\alpha\in\srts_Q\mid w^{-1}\alpha<0\}}$ times the
characteristic function of the projection to $\aaa_Q$ of the set
\[
T+\{\sum_{\alpha\in\srts_0}c_\alpha\alpha^\vee\mid c_\alpha>0\iff w^{-1}\alpha<0\}.
\]
In particular, if $\psi\in\AF_G^2$ and $P\ne G$, then
\begin{equation} \label{eq: ortheis}
\inner{\Lambda^TE(\varphi,\lambda)}{\psi}_{\autspace_G}^{\unit}\equiv0.
\end{equation}
\end{proposition}

The left-hand side of \eqref{eq: main34} is an analytic function on $\iii\aaa_P^*$ (viewed as a distribution).
On the right-hand side of \eqref{eq: main34}
\[
\opr:C^\infty(\iii\aaa_Q^*)\otimes\Polyexp^{\unit}(\aaa_Q)\rightarrow\operatorname{End}(\mathcal{D}(\iii\aaa_Q^*))
\]
is the linear map (rather than homomorphism) with $C^\infty(\iii\aaa_Q^*)$ acting by multiplication and $p\in\Polyexp^{\unit}(\aaa_Q)$
acting by convolution by $\hat p$.

\begin{remark}
Taking $\psi$ to be an Eisenstein series, we obtain Arthur's asymptotic inner product of truncated
Eisenstein series (cf.\ \cites{MR650368, MR2767521}).\footnote{Arthur's original proof
used Langlands's description of the discrete spectrum.}
\end{remark}

The first part of Theorem \ref{thm: main2} immediately follows from Proposition \ref{prop: maincor}
by taking the limit in $T$ as $\min_{\alpha\in\srts_0}\sprod{\alpha}T\rightarrow\infty$,
since the limit (as tempered distributions) of $\chi_w^T$ is the constant function $1$ if $w=e$ (and $Q=P$) and $0$ otherwise.
For the second part, if $\psi=E(\varphi',\mu)$ then
\[
\psi_P^{\weak}=\sum_{w\in\Weyl^{\supset Q}(P)}E^P(M(w^{-1},\mu)\varphi',w^{-1}\mu),
\]
and the statement follows from \eqref{eq: ortheis} by taking $T\rightarrow\infty$.

\subsection{}

The proof of Proposition \ref{prop: maincor} is based on the following generalization of
a result of Arthur on the truncation of Eisenstein series induced from cuspidal automorphic forms \cite{MR558260}*{Lemma 4.1}.\footnote{This formula,
in the cuspidal case, was in fact Langlands's original definition of truncated Eisenstein series. Subsequently, Arthur realized that
a truncation operator can be defined more generally.}

\begin{proposition}[\cite{MR3001807}] \label{prop: trn12}
For any $\varphi\in\AF_P$ and $\Rel\lambda\gg0$, $\Lambda^TE(\varphi,\lambda)$ is the sum over $Q\in\stds$ and
$w\in\Weyl^{\supset Q}(P)$ of
\[
\sum_{\gamma\in Q(F)\bs G(F)}
\Lambda^{T,Q}(M(w,\lambda)\varphi_{P_w})_{w\lambda}(\gamma g)\chi_w(H_Q(\gamma g)-T).
\]
\end{proposition}
Here, $\Lambda^{T,Q}$ is a relative truncation operator and $\chi_w$ is the function on $\aaa_Q$ given by
\[
\newchar_w(X)=\begin{cases}(-1)^{\abs{\posrts QX}}&
\begin{aligned}\text{if \eqref{eq: prop22} holds},
\end{aligned}
\\0&\text{otherwise}\end{cases}
\]
where $\posrts Q{\sum_{\alpha\in\srts_Q}x_\alpha\alpha^\vee}=\{\alpha\in\srts_Q\mid x_\alpha>0\}$
and the condition \eqref{eq: prop22} is
\begin{equation}\tag{*} \label{eq: prop22}
\posrts QX=\{\alpha\in\srts_Q\mid w^{-1}\alpha<0\text{ or }w^{-1}\alpha\in\srts_{P_w}^P\text{ and }\sprod{\alpha}X\le0\}.
\end{equation}

In particular, $\newchar_w(\cdot-T)=\chi_w^T$ if $w\in\Weyl(P,Q)$.
In general, $\newchar_w$ can be expressed as a linear combination of characteristic functions of simplicial cones.
Thus, its Laplace transform is a rational function with hyperplane singularities.

Using Proposition \ref{prop: trn12} and simple geometric properties of $\chi_w$, we can compute the inner product
\[
\inner{\Lambda^TE(\varphi,\lambda)}{\psi}_{\autspace_G}
\]
as a sum of contributions over $Q\in\stds$ and $w\in\Weyl^{\supset Q}(P)$. One sees that only
$w\in\Weyl(P,Q)$ may contribute to
\[
\inner{\Lambda^TE(\varphi,\lambda)}{\psi}_{\autspace_G}^{\unit}.
\]
This contribution can be computed explicitly, and yields Proposition \ref{prop: maincor}.

\begin{remark} \label{rem: eisholom}
Although we gave ourselves the liberty of using the analyticity of the Eisenstein series on $\Rel\lambda=0$,
one can run the argument without this assumption. Namely, Proposition \ref{prop: maincor} and Theorem \ref{thm: main2} hold
provided that $f$ vanishes to a high order on the hyperplane singularities of the Eisenstein series.
Then, by an easy argument one concludes the analyticity of the Eisenstein series on $\Rel\lambda=0$ \cite{MR2767521}.
\end{remark}

\section{Fourier coefficients and restrictions of Eisenstein series}

One of the precursors of the modern theory of automorphic forms is the classical Siegel--Weil formula,
which relates
Fourier coefficients of Eisenstein series on the symplectic group induced from
the trivial character of the Siegel parabolic to representation numbers of quadratic forms
(see e.g.\ \cite{MR2402685}).
This Eisenstein series is also the linchpin for the doubling method (see below).
Its derivative plays a key role in Kudla's program \cite{MR2083214}.
While degenerate Eisenstein series have other striking applications as well
(e.g.\ \cite{MR2341906}) we will mostly focus here on Eisenstein series induced from cuspidal representations.

\subsection{Constant term}
In his famous Yale manuscript \cite{MR0419366}, Langlands computed the constant term of an Eisenstein series
induced from a cusp form on a maximal Levi subgroup $M$ of $G$.
This led him to the concept of the $L$-group, the notion of a general automorphic $L$-function
and ultimately, to the principle of functoriality and the birth of the Langlands program \cite{MR0302614}.
The importance of this computation cannot be overstated.

Langlands's computation, together with the meromorphic continuation of Eisenstein series, implies the meromorphic continuation of
automorphic $L$-functions $L(s,\pi,r_i)$, $i=1,\dots,m$ pertaining to (irreducible) cuspidal representations $\pi$
of $M(\A)$ and the irreducible constituents $r_1,\dots,r_m$ of the adjoint representation of ${}^LM$ on the Lie algebra
of the unipotent radical of ${}^LP$. These include Rankin-Selberg $L$-function
for $\GL_n\times\GL_m$, tensor product $L$-function for $G\times\GL_n$ where $G$ is a classical group
(with respect to the canonical representation of ${}^LG$), symmetric square, exterior square
and twisted tensor $L$-functions of $\GL_n$ and a finite exotic list of cases pertaining to exceptional groups.
Moreover, since the Eisenstein series are meromorphic functions of finite order, the same holds for the $L$-functions above \cite{MR2230919}.

On the other hand, Langlands conjectured not only that a general automorphic $L$-function can be meromorphically continued
to $\C$, but also that it admits a functional equation and only finitely many poles
(and in fact, is equal to the standard $L$-function of an automorphic form on $\GL_n$).
These properties seem out of reach of the constant term computation, even for the $L$-functions that fall under its jurisdiction.

\subsection{Nondegenerate Fourier coefficients}
Assume now that $G$ admits a maximal unipotent subgroup $U_0$ (namely, the radical of a minimal parabolic subgroup of $G$)
that is defined over $F$, i.e., that $G$ is quasi-split.
Shahidi computed the non-degenerate Fourier coefficients with respect to $U_0$ of an Eisenstein series as above.
This is zero unless $\pi$ itself is generic, i.e., admits a non-zero non-degenerate Fourier coefficient
with respect to $U_0\cap M$. Assuming from now on that $\pi$ is generic, and using the Casselman--Shalika formula,
one recovers (up to local factors) the denominator in Langlands's formula \cite{MR610479}.
In the case of the Borel subgroup, this computation goes back in fact to Jacquet's thesis \cite{MR0271275}.

The Langlands--Shahidi method has far-reaching applications, thanks to the work of Cogdell, Kim, Piatetski-Shapiro
and of course Shahidi himself.
We refer the reader to \cite{MR2683009} and the references therein for a complete account, including many local aspects
(or \cite{MR3307916} for a more concise account) and confine ourselves to a very brief discussion.
First, it provides a functional equation for the $L$-functions $L(s,\pi,r_i)$ above.
Second, it implies that the poles (resp., zeros) of the partial $L$-functions for $\Rel s\ge\frac12$
(resp., $\Rel s\ge1$) are all on the real line, and in particular there are only finitely many of them in these regions.
Also, in the region $\Rel s\ge\frac12$ apart from the finitely many poles, the partial $L$-functions are bounded polynomially.
With considerable more input from the local theory \cite{MR3068402}, one shows that in fact the completed $L$-functions have only finitely many poles.
Moreover, if $(M,\pi)$ does not admit a non-trivial symmetry under conjugation by $G$ (which can always be arranged upon twisting $\pi$
by a sufficiently ramified character) then the Eisenstein series is holomorphic for $\Rel s\ge0$ and the completed
$L$-functions are entire. This (together with the aforementioned polynomial growth) is crucial for the application of the converse theorem.
The upshot is a proof of functorial transfer of generic representations of classical groups to $\GL_n$,
as well as the functorial transfer from $\GL_2\times\GL_3$ to $\GL_6$. (The latter yields the best known bounds towards
the Ramanujan--Selberg conjectures for Maass forms.)

We note that although Langlands's constant term computation carries over to the case of covering groups \cite{MR3749191},
non-degenerate Fourier coefficients of Eisenstein series on covering groups exhibit distinct new features.
They were studied extensively by Brubaker, Bump, Friedberg and others. In particular, they discovered
unexpected and deep connections to physics.
We refer the reader to \cite{MR3385631}, the survey \cite{MR2896085} and the references therein.

\subsection{Other constructions}
The constant term and non-degenerate Fourier coefficients of Eisenstein series are only the tip of the iceberg.
The Rankin--Selberg method provides invaluable information about $L$-functions using other constructions involving Eisenstein series.\footnote{In
his overview \cite{MR1476510}, Langlands was highly critical of the Rankin--Selberg method for its perceived lack of structure.
Perhaps the situation is changing now through the works of Yiannis Sakellaridis and others \cite{2111.03004}.}

For instance, let $G$ be a classical group, i.e., the isometry group of $(V,q)$ where $q$ is a symmetric/antisymmetric/hermitian form on
a finite-dimensional vector space $V$. The isometry group $H$ of $(V\oplus V,q\oplus -q)$ is quasi-split and we may consider
the Eisenstein series $E_\chi(s)$ on $H$ induced from a one-dimensional character $\chi\circ\det$ of the Siegel parabolic subgroup of $H$ (whose Levi subgroup
is $\GL(V)$). The cuspidal part of the spectral expansion of the restriction of $E_\chi(s)$ to $(G\times G)(F)\bs(G\times G)(\A)$ is supported on the direct sum of $\pi\otimes\pi^\vee$
where $\pi$ ranges over the cuspidal automorphic representation of $G(\A)$, and the coefficients involve the values $L(s+\frac12,\pi\otimes\chi)$
of $L$-functions corresponding to the canonical representation of ${}^LG\times\mathbb{G}_m$.
This is the contents of the doubling method of Piatetski-Shapiro and Rallis \cites{MR892097, MR2192828}.

Another striking example is the descent construction of Ginzburg--Rallis--Soudry \cite{MR2848523}.
To fix ideas, consider a cuspidal representation $\pi$ of $G'=\GL_{2n}$ such that the exterior square $L$-function has a pole at $s=1$.
Then, the Eisenstein series on the orthogonal group $\SO(4n)$ induced from $\pi\abs{\det}^s$ on the Siegel parabolic has a pole at $s=\frac12$.
Taking a certain Fourier coefficient, stabilized by $G=\SO(2n+1)$, of the residue representation, one obtains an automorphic representation $\sigma$ of $G(\A)$.
The main result is that $\sigma$ is the generic cuspidal representation of $G$ whose functorial transfer to $G'$ is $\pi$.
Among other things, this leads to a relation between the Petersson inner product and the non-degenerate Fourier
coefficient of generic cusp forms \cite{MR3649366}.
More recently, Ginzburg and Soudry studied a certain Fourier coefficient, stabilized by $G\times G$, of an iterated residue of an Eisenstein series on $\SO(4n(2n+1))$
induced from $\pi\otimes\dots\otimes\pi$ ($2n+1$ times) \cite{MR4389019}.
They obtain the direct sum of $\sigma\otimes\sigma^\vee$ where $\sigma$ ranges over the entire $L$-packet indexed by $\pi$.
This gives an altogether different approach to Arthur's endoscopic classification of representations of classical groups
without appeal to the trace formula. It may also lead to a new strategy for the proof of the Gross--Prasad conjectures
in quantitative form, including in cases which so far resisted the relative trace formula approach.

The study of Fourier coefficients of Eisenstein series and their residues is extensive but at the same time mysterious \cite{MR2214128}.
So far it has been mostly driven by ``the proof of the pudding is in the eating'' approach, with unquestionable success.
It would be desirable to perceive a more conceptual perspective, at least to be able to predict when such coefficients are potentially meaningful.

I will not discuss here regularized integrals of Eisenstein series over reductive subgroups.
This is technically important for the relative trace formula and other applications.
Recently, it was developed in great generality by Zydor \cite{1903.01697}.
Likewise, the computation of such integrals for residues of Eisenstein series is a very fruitful and well-studied theme
with an abundance of applications. We refer the reader to \cite{MR4255059} and the references therein for some new and old results in this direction.

\section{Open questions}

Let me end this report by mentioning some vaguely formulated (mostly folklore) open problems and questions.

\subsection{Poles of Eisenstein series}

Suppose that $\pi$ is a cuspidal representation of the Levi subgroup of a maximal parabolic subgroup of $G$.
What can be said about the poles of the Eisenstein series induced from $\pi$ for $\Rel s>0$ ?\footnote{The poles
for $\Rel s<0$ are related to zeroes of $L$-functions, and thus to the generalized Riemann Hypothesis.
See \cites{MR633666, MR634284} for a relation (via Eisenstein series) between RH and the rate of equidistribution of horocyles in $\SL_2(\Z)\bs\HHH$.}

By general theory, it is known that these poles are simple and real and they can only occur
if the data $(M,\pi)$ is invariant under the long Weyl element \cite{MR1361168}*{Ch.\ 4}.
In particular, there are only finitely many of them.
It is expected that the poles\footnote{with the usual normalization where $1$ is identified with the fundamental weight.} are contained in $\frac12\mathbb{N}$.
In particular, their number is bounded in terms of $G$ only.\footnote{In the non-arithmetic case, the number of poles
cannot be bounded independently of the lattice.}
This is known at least for the general linear group \cite{MR1026752} and (using Arthur's work) for classical groups \cite{MR2457189}*{Proposition 1.2.1}.
However, this is open in general, even in the case where $\pi$ is generic, in which it is expected that the only possible pole is at either $\frac12$ or $1$.
(Even less is known if $\pi$ is only assumed to be in the discrete spectrum \cite{MR1720186}.)
In any case, the situation is somewhat unsatisfactory.
Current methods are based on a detailed analysis of $L$-functions and normalized intertwining operators.
It would be desirable to have a ``soft'' general spectral method to control the number of poles (if not their location).

A related problem is a geometric normalization of Eisenstein series.
In the case of Eisenstein series on $\GL_{n+m}$ induced from a cuspidal representation of $\GL_n$
and a character of $\GL_m$, one can normalize the Eisenstein series using a family of holomorphic sections
obtained from Schwartz function on an affine space (see \cite{MR0432596} for the case $m=1$). No such procedure is known is general. (See \cites{MR1694894, MR1988971}
for some special cases.)

\subsection{The discrete spectrum}
The residues discussed above furnish the contribution to the discrete spectrum from Eisenstein series in corank one.
The contribution of other Eisenstein series to the discrete spectrum is much more complicated.

What can be said about them in general?

There is an extensive work by M\oe glin \cites{MR1091891, MR1159268, MR2457189, MR2822218} and others on classical groups.
More recent works by de Martino--Heiermann--Opdam \cite{1512.08566} and Kazhdan--Okounkov \cite{2203.03486}
give some hope for a more uniform
and conceptual approach, but so far they are restricted to the unramified case.
In any case, it is fair to say that the situation is far from understood.

A basic open question in this regard is the validity of the upper bound in Weyl's law.
Thus, fixing the level, we would like to know that the counting function for the discrete spectrum of the locally symmetric space
with Laplace eigenvalue $\le T$ obeys the upper bound of Weyl's law as $T\rightarrow\infty$.
In fact, for the cuspidal spectrum we already know the Weyl law, even with a power saving remainder term
\cites{MR2306657, MR4280092}, so the ``only'' issue is
to show that the noncuspidal discrete spectrum is negligible. Embarrassingly, this is still unknown in general.

\subsection{Paley--Wiener Theorems}
Plancherel formulae are usually proved by performing analysis on a suitable space of test functions
(such as the Schwartz space), rather than on the $L^2$-space itself.
They come hand-in-hand with corresponding results on these spaces of functions.
This is also true for Langlands's spectral decomposition of $L^2(G(F)\bs G(\A))$
(cf.\ \cite{MR1603257}*{Theorem 12}).
However, the result in [loc.\ cit.] is for the space of $K$-finite, spectrally bounded functions.
It would be advantageous to consider larger spaces, such as the Harish-Chandra Schwartz space.
This problem is wide open. (See \cite{MR3156857} for a discussion and some partial results for the general linear group.)
An analogous problem for the (much smaller) space of smooth functions of rapid decay is discussed by Casselman. (See \cite{MR2201004} and the references therein.)

A closely related problem is the growth of Eisenstein series, a problem considered most recently by
Subhajit Jana and Amitay Kamber.
A key issue is understanding the inner product of truncated Eisenstein series.
They are expressed as a sum over the Weyl group elements of certain divided differences.
It is possible to evaluate the limits of the summands
in terms of rank one intertwining operators and their first-order derivatives \cites{MR2811598, MR2811597}.
This is sufficient for applications to the trace formula.
However, in order to understand the truncated inner product in general, it is necessary to
control all values, not just the limiting values. (In the noncuspidal case there are additional difficulties.)

\subsection{Lower bounds on $L$-functions}
It has been known for some time that the holomorphy of Eisenstein series on the unitary axis is closely related
to non-vanishing of $L$-functions on $\Rel s=1$ \cite{MR0432596}.
In \cite{MR2058625} Sarnak made this relation quantitative.
In particular, the classical Eisenstein series \eqref{eq: ezs} is used to obtain a zero-free region for the Riemann zeta function
that is roughly on par with the classical one due to de la Vall\'ee Poussin.
The core of Sarnak's idea is a use of Bessel's inequality to bound sums of squares of Fourier coefficients of $E(z;s)$
by the inner product of $E(z;s)$ (suitably truncated) which in turn is computed by the Maass--Selberg relations.
The inhomogeneous feature of the ensuing inequality yields a lower bound on $\abs{\zeta(1+\iii t)}$ in terms of an upper bound
on $\zeta'(1+\iii t)$ (which is easy). (In order to optimize the bound, one also needs a sieve-theoretic argument.)
The argument, at least in its crude form, extends to $L$-functions of generic representations
of Levi subgroups for which the Langlands--Shahidi method applies \cite{MR2230919}.
There is plenty of room to sharpen these bounds and to extend this method beyond generic representations.
We point out that for Rankin--Selberg $L$-functions this method seems inferior to classical methods \cite{MR3980284}.

This is just a sample of problems and questions. It is clear that Eisenstein series will continue to be a mainstay
of the theory of automorphic forms for the foreseeable future.

%\bibliographystyle{amsplain}
%\bibliography{../Bibfiles/all}

\def\cprime{$'$}

\end{document}